\documentclass{article}

\usepackage{amsmath}
\usepackage{amssymb}

\newtheorem{Coro}{Corollary}

\newtheorem{Lemm}{Lemma}
\newtheorem{Prop}{Proposition}
\newtheorem{Theo}{Theorem}

\begin{document}
\centerline{\bf\sc\Large Logarithmic Trace of Toeplitz Projectors}

\medskip\hskip55mm L. Boutet de Monvel

\begin{abstract}In \cite{BG81} we defined Toeplitz projectors on a compact contact manifold, which are analogues of the Szeg\"o projector on a strictly pseudo-convex boundary. The kernel of a Toeplitz projector, as the Szeg\"o kernel, has a holonomic singularity including a logarithmic term. The coefficient of the logarithmic term is well defined, so as its trace (the integral over the diagonal). Here we show that this trace only depends on the contact structure and not on the choice of the Toeplitz operator (for a given contact structure there are many possible choices). This generalizes a result of K. Hirachi \cite{Hi04} for the Szeg\"o kernel, and also shows that his invariant (the trace of the logarithmic coefficient of the Szeg\"o kernel) only depends on the contact structure defined by the  boundary pseudo-convex CR structure. Finally we show that the Toeplitz logarithmic trace vanishes identically for all contact forms on the three-sphere.\end{abstract}

%%%%%%%%%%%  1  %%%%%%%%%%%%%
\section{Introduction}

Let $(X,\lambda)$ be a compact oriented contact manifold of dimension $2n-1$. This means that $X$ is a manifold equipped with a differential 1-form $\lambda$ such that $\lambda(d\lambda)^{n-1}$ is a volume element  ($2n-1={\rm dim\,}X$); two forms $\lambda,\lambda'$ define the same structure iff $\lambda'=f\lambda$ with $f$ a smooth positive function. Equivalently $X$ is equipped with a smooth symplectic half-line sub-bundle $\Sigma\subset T^*X$ (the set of positive multiples of $\lambda$).

A typical example is the unit sphere $X\ (z.\bar{z}=1)$ in $\mathbb{C}^n$, with $\lambda=\frac1i\bar{z}.dz|_X$. The Szeg\"o kernel of the unit sphere is 
$$
S=\frac1c(1-z.\bar{w})^{-n}\qquad \mbox{with}\quad c=\frac{2\pi^n}{(n-1)!}, \quad\mbox{the volume of }X
$$
It is the kernel of the orthogonal projector (still denoted $S$) on the space of holomorphic functions (${\rm ker\,}\bar{\partial}_b$) in $L^2(X)$:
$$
Sf(z) = \int_X S(z,\bar{w}) f(w) d\sigma(w),\qquad \mbox{with } d\sigma(w)\mbox{ the canonical volume element of }X.
$$
The Szeg\"o kernel is linked to the contact structure in the following manner: $S$ is a Fourier integral operator with complex canonical relation $\mathcal{C}$, where the real part of $\mathcal{C}$ is ${\rm Id}_\Sigma$, the graph of the identity map of $\Sigma$ (which just means that in the complexification of $X\times X$ the set of real points of the hypersurface $\{\phi=i(z.\bar{w}-1)=0\}$ is the diagonal $z=w$, and there we have $d_z\phi =d_{\bar{w}}\phi=c\lambda$, with $c>0$ - in this case $c=1$).

\medskip Things are similar when $X$ is the boundary of a complex manifold $\Omega$ and is strictly pseudo-convex. In this case $S$ is again defined as the orthogonal projector on the space of boundary values of  holomorphic functions (${\rm ker\,}\bar{\partial}_b$) in $L^2(X,dv)$ (this depends on the choice of a volume element $dv$). If $X$ is defined by a real analytic equation $q(z,\bar{z})=0$ ($q>0$ in $\Omega$), the Szeg\"o kernel $S$  is holomorphic w.r. to $z$ and antiholomorphic w.r. to $w$, smooth up to the boundary except along the diagonal of $X\times X$ where it has a typically holonomic singularity (cf. \cite{lB76.5}, \cite{mK77}):
\begin{equation} \label{T2}
S=\phi(z,\bar{w})\,(q(z,\bar{w})+0)^{-n}+\,\psi(z,\bar{w})\,{\rm Log\,}(q(z,\bar{w})+0)
\end{equation}
where $\phi$ and $\psi$ are smooth functions defined near the diagonal, holomorphic (resp. anti--) w.r. to $z$ (resp. to $w$). If $X$ is only $C^\infty$ rather than real analytic, the formula above is formal, i.e. it only defines  the Taylor series of $q,\phi$ and $\psi$ along the diagonal, but this is enough to determine completely the singularity of $S$.

\bigskip
The coefficient $\psi$ of the logarithmic term (or rather its Taylor series along the diagonal) is completely defined; it only depends on the complex boundary (CR) structure of $X$ and $dv$. In \cite{Hi04} Hirachi shows that its trace, i.e. the integral
\begin{equation}\label{Ltr}
L(\lambda) = \int_X \psi(z,\bar{z}) \ dv
\end{equation}
is a rather rigid invariant: it only depends on the $CR$ structure of $X$, not on the choice of the volume element $dv$, and it is also invariant under deformation.

\bigskip We have shown in \cite{BG81} that if $X$ is an oriented contact manifold, there always exists a
Toeplitz projector, which is an analogue of the Szeg\"o projector: this is  a projector in $L^2(X,dv)$ (for some volume element $dv$) which is an ``elliptic" Fourier integral operator with complex canonical relation $\mathcal{C}$, where again the real part of $\mathcal{C}$ is the graph of ${\rm Id}_\Sigma$. As in the holomorphic situation, $\mathcal{C}$ is the conormal bundle of a complex hypersurface $\{q(x,y)=0\}$ in the complexification of $X\times X$, $q$ a smooth function, with $q=0$, $id_xq=-id_yq=c\lambda$, a positive multiple of the contact form, on the diagonal, and ${\rm Re\,}q\ge c |x-y|^2$ near the diagonal ($c>0$; ${\rm Re\,}q>0$ outside of the diagonal; there are further conditions to ensure $\mathcal{C}\circ\mathcal{C}=\mathcal{C}$). 

The kernel of $S$ is a Fourier integral:
\begin{equation} \label{T1}
S(x,y) \sim \int_0^\infty e^{-Tq(x,y)} a(x,y,T)\ dT
\end{equation}
with $a\sim\sum_{k<n} a_k(x,y) T^k$, a symbol of degree $n-1$. $S$ has a holonomic singularity as in \eqref{T2}, which we can expand as
\begin{equation} \label{T3}
S(x,y)\sim \sum_{0<k\le n}\alpha_k(x,y)(q(x,y)+0)^{-k}+\sum_{k\ge0}\beta_k(x,y) q(x,y)^k{\rm Log\,}(q(x,y)+0) \quad({\rm mod. }C^\infty)
\end{equation}
with 
%\alpha_k(x,y)=\frac{a_{k-1}(x,y)}{(k-1)!}\,\quad{\rm for}\ k>0, \quad \beta_k(x,y)= \frac{a_{-1-k}(x,y)} {k!} \quad{\rm for}\ k\ge0
\begin{equation} \label{T4}
\alpha_k(x,y)=a_{k-1}(x,y)\,(k-1)!^{-1} \quad{\rm for}\ k>0, \quad \beta_k(x,y)=a_{-1-k}(x,y)\,k!^{-1} \quad{\rm for}\ k\ge0
\end{equation}
Here again the coefficient $\psi(x,y)$ of the logarithmic term is well defined; it depends on the choice of the projector $S$. In this article we prove

\begin{Theo}\label{th1}  The trace $L(S)=\int_X\psi(x,x) dv$ of the logarithmic term only depends on the contact structure of $X$, not depend on the choice of the of the canonical relation $\mathcal{C}$ of $S$.\end{Theo}

\medskip To prove the theorem we will use the fact that any two Toeplitz projector singularities belonging to the same contact structure can deformed one into the other. Note that the theorem also implies that the trace of the logarithmic coefficient is invariant in a deformation of the contact structure, because such deformations are always trivial: if $\lambda_t$ is a smooth one-parameter family of contact structures on $X$ ($X$ compact), there exists a smooth family $\Phi_t$ of diffeomorphisms of $X$ such that $\Phi_t^*\lambda_t$ is a multiple of $\lambda_0$, and of course diffeomorphisms preserve separately polar and logarithmic singularities.

The result of K. Hirachi for the Szeg\"o kernel follows from theorem \ref{th1}, which shows moreover that $L(\lambda)$ only depends on the contact structure. I have no example of a contact structure $(X,\lambda)$ with $L(\lambda)\neq0$ and Hirachi's question on that point remains open, but  in section \ref{sph}, I show that if $X$ is the 3-sphere, $L(\lambda)$ always vanishes. In view of this it not unlikely that $L(\lambda)$ vanishes for all contact forms.

%%%%%%%%%%%  2  %%%%%%%%%%%%%
\section{Toeplitz projectors}\label{deform}

As mentioned above a Toeplitz projector  is a projector $S$ in $L^2(X,dv)$ which is an elliptic Fourier integral operator with positive complex canonical relation $\mathcal{C}$ such that  real part is the is the graph of $Id_\Sigma$.  It can be represented as a Fourier integral \eqref{T1} with , $q$ a smooth function on $X\times X$, ${\rm Re\,}q\ge c\;|x-y|^2$ near the diagonal, and $id_xq=-id_yq=c\,\lambda$, ($\lambda$ the contact form; there are further conditions on $q$ ensuring $\mathcal{C}\circ\mathcal{C}=\mathcal{C}$). The construction of such canonical relations and Toeplitz projectors is described in  \cite{BG81}, see also \cite{lB86.2},  \cite{lB97.3}. This construction allows deformations or compact group actions and in fact shows that set of such canonical relations, so as the set of singularities of Toeplitz projectors, is contractible; in particular we have:

\begin{Prop} Let $\mathcal{C}_0,\mathcal{C}_1$, resp. $S_0,S_1$ be two canonical relations, resp.  Toeplitz projectors as above. Then there exist smooth one parameter families  $\mathcal{C}'_t$, resp. $S'_t$ of canonical relations, resp. Toeplitz projectors such that $\mathcal{C}'_i=\mathcal{C}_i$, resp. $S'_i-S_i$ has a smooth kernel, for $i=0,1$.
\end{Prop}

Another way of stating this is
\begin{Coro} \label{soft}
Let $(X,\lambda)$ be an oriented contact manifold. Then the singularities of Toeplitz kernels form a soft sheaf on diag$\,X\subset X\times X$.
\end{Coro}
In other words if $S_0$ is a Toeplitz kernel singularity defined in a neighborhood of some closed subset $K\subset X$, there exists a global singularity $S$ (defined on the whole of $X$) which coincides with $S_0$ near $K$. {\small(The mod. $C^\infty$ assertion about $S$ cannot be improved: the set of Toeplitz projectors is not connected; in fact if $S,S'$ are two Toeplitz projector, $S\circ S'$ induces a Fredholm operator from the range of $S'$ to the range of $S$, whose index is an arbitrary integer and is of course deformation invariant.
The assertions above are possibly easier to see in the setting of "Hermite operators" of \cite{lB74.4} or of "symplectic spinors" of \cite{vG74}: the leading term is parametrized by normalized gaussian densities of the form $({\rm discr}\,2\pi q)^{-\frac1/2}\exp-\frac12q$ with $q$ a quadratic form, ${\rm Re\,}q\gg0$, and the rest is a formal expansion easily linearized. But in this setting it is more awkward to keep track of the logarithmic term.)}

\medskip Since these facts are not explicitly stated in loc. cit., I briefly recall the proof in \S\ref{Toepl} below.

%%%%%%%%%%%  3  %%%%%%%%%%%%%
\section{Logarithmic term of the Toeplitz projector}

Let $S$ be a Toeplitz projector, defined by a Fourier integral as in \eqref{T1}
$$
S(x,y) \sim \int_0^\infty e^{-Tq(x,y)} a(x,y,T)\ dT
$$
Once $q$ and $a$ are fixed, following Hirachi, we introduce the asymptotic expansions
\begin{equation}\label{TA1}
S(x,y,\epsilon) = \int_0^\infty e^{-T(q+\epsilon)}a\,dT\sim\sum_{k>0} \alpha_k(x,y)(q+\epsilon)^{-k}+\sum_{k\ge0}\beta_k(x,y)(q+\epsilon)^k{\rm Log\,}(q+\epsilon)
\end{equation}
and its trace 
\begin{equation}\label{TA2}
s(\epsilon)={\rm tr\,}S_\epsilon=\int_0^\infty e^{-T\epsilon} a(x,x,T) dv(x)dT\sim \sum \alpha_k\epsilon^{-k}+\sum\beta_k\epsilon^k{\rm Log\,}\epsilon
\end{equation}
with $\alpha_k(x,y),\beta_k(x,y),\beta_k$ as in \eqref{T4}, $ \alpha_k =\int  \alpha_k(x,x), \beta_k=\int \beta_k(x,x)$.

\medskip These should be understood as a Fourier integral distribution on $X\times X\times \mathbb{R}$ defining an $\epsilon$-parameter family of operators on $X$, resp. on $\mathbb{R}$, with complex Lagrangian the conormal bundle $T^*_{\{q+\epsilon=0\}}(X\times X\times \mathbb{R})$, resp. $T^*_0\mathbb{R}$. Although as distributions they are only defined for $\epsilon\ge 0$, they satisfy regular holonomic systems of pseudo-differential equations, with well defined jets of infinite order along ${\rm diag\,}X\times \{\epsilon=0\}$, (resp. at $\epsilon=0$); this is the only thing that counts. 

\medskip The logarithmic coefficients $\beta_0(x,y),\beta_0$ only depend on $S$, but of course the full asymptotic expansions $S_\epsilon(x,y), s(\epsilon)$ depend on the precise choice of the phase function $q$ and symbol $a$ (the leading term transforms in an obvious way if we replace $q$ by a multiple of $q$, but the rest depends in a more complicated manner on the choice of $\phi,a$). For the proof it will be convenient to choose $\phi$ and $a$ adequately:

\begin{Lemm} \label{Ldef} 
We can choose the phase $q$ and symbol $a$ so that $S_\epsilon\sim S\circ S_\epsilon \circ S$. \end{Lemm}
Proof: the theorem of the stationary phase shows that we have 
$$
S\circ e^{-\lambda q} a \circ C \sim \int e^{-Tq(x,u)}a(x,u,T)e^{-\lambda q(u,v)}a(u,v,\lambda) e^{-xq(v,y)}a(v,y,S) dTdSdv(u)dv(v)
$$
$$
\sim e^{-\lambda q'(x,y)} a'(x,y,\lambda)
$$
for some phase function $q'(x,y)$ and symbol $a'(s,y,\lambda)$ ($\lambda q'$ the critical value of $-Tq(x,u)_\lambda q(u,v)-Sq(v,y)$ at its stationary point, $a'$ given by the asymptotic expansion in the method of the stationary phase). The projector equation $S\circ S = S$ then implies that
$$
S\sim \int e^{-\lambda q'(x,y)} a'(x,y,\lambda)d\lambda,
$$
and 
$$
S\circ e^{-\lambda q'(x,y)} a'(x,y,\lambda)\circ \sim e^{-\lambda q'(x,y)} a'(x,y,\lambda)
$$
In particular $q'$ is a smooth multiple of $q$. In the holomorphic case (Szeg\"o kernel) this just means that $q'$ and $a'$ are holomorphic in $x$ and antiholomorphic in $y$.

\begin{Prop} 
Let $S_t$ be a smooth family (deformation) of Toeplitz projectors. Then the Logarithmic coefficient $\beta_0(t)$ is constant.
\end{Prop}
Proof: we write $S$ for $S_t,D=\frac d{dt}$. We have $Ds(\epsilon)=\sum D \alpha_k\epsilon^{-k} +\sum D\beta_k\epsilon^k {\rm Log\,}\epsilon$.
Since $S$ is a projector, we have $DS =DS\;S+S\;DS$, hence
\begin{equation} \label{derivS}
DS=[S,(2S-1)DS] .
\end{equation}
Let us choose choose the smooth phase function $q_t$ and symbol $a_t$ as in Lemma\ref{Ldef}, so that $S S_\epsilon \sim S_\epsilon \sim S_\epsilon S$. With this choice, the kernel of $[S,(2S-1)DS_\epsilon]$ is a Fourier integral distribution belonging to the same Lagrangian as $S_\epsilon$ (the conormal bundle of the hyper-surface $\{q_t+\epsilon=0\}$) and $DS_\epsilon - [S,(2S-1)DS_\epsilon]$ vanishes for $\epsilon=0$ , i.e. it is a multiple of $\epsilon$, a Fourier integral of the form
$$
 \int e^{-Tq(x,y)+\epsilon}\ \epsilon\ b(x,y,\epsilon,T)\;dT
$$
Since the trace of a commutator vanishes identically, we have 
$$
Ds(\epsilon)\sim \int e^{-Tq(x,y)+\epsilon}\ \epsilon\; b(x,y,\epsilon,T)\;dT\,dv.
$$
Now we can repeat the argument of K. Hirachi \cite{Hi04}: the asymptotic expansion of $s(\epsilon)$  is a multiple of $\epsilon$, so in it the coefficient of ${\rm Log\,}\epsilon$ vanishes. Since any two Toeplitz projectors can be deformed into one another mod. smoothing operators, this proves theorem \ref{th1}.

%%%%%%%%%%%  4  %%%%%%%%%%%%%
\section{Deformations of Toeplitz projectors.}\label{Toepl}

As recalled above, a Toeplitz projector is a Fourier integral operator with complex phase function. Its kernel is of the form \eqref{T1}:
$$
S(x,y) = \int_0^\infty e^{-Tq(x,y)} a(T,x,y)dT = \phi(x,y) (q+0)^{-n}+\psi(x,y){\rm Log\,}(q+0)
$$
where $a$ is a symbol of degree $n-1$, $\phi,\psi$ are smooth functions. The phase function is $\Phi=iTq$, with $q$ a smooth complex function on $X\times W$, $q(x,x)=0,dq\neq0$ on the diagonal, ${\rm Re\,}q(x,y)\ge c|x-y|^2$  near the diagonal.

The objects we are dealing with are really jets of infinite order along $\Sigma$ in $T^*X$, the diagonal in $X\times X$ or the diagonal of $\Sigma$ in $T^*(X\times X)$, but it will be be more agreeable, and perfectly legitimate, to use the language of functions and sub-manifolds of differential geometry.  

The complex canonical relation $\mathcal{C}$ is the set of covectors $(x,\xi)=d_x\Phi$, $(y,\eta)=-d_y\Phi$ with $q(x,y)=0$, in the complexification of $T^*X\times T^*X^0$. Its real part is ${\rm diag\,}\Sigma$, i.e. $\frac1id_xq=-\frac1id_yq=c\lambda$ on the diagonal ($\lambda$ the contact form, $c>0$), with $c>0$. The positivity of $\mathcal{C}$ follows from the condition ${\rm Re\,}q\ge c|x-y|^2$ (cf. \cite{MS74}).  For a given projector $S$ only the complex canonical relation, i.e. the hypersurface $\{q=0\}$ is well defined; the phase or $q$ itself is only defined up to a smooth factor.  We will not require that $q$ be hermitian ($q(y,x)=\bar{q}(x,y)$) or $S$ self adjoint - this anyway makes sense only once the volume $dv$ is chosen.

%%%%%%%%%%%  4.1  %%%%%%%%%%%%%
\subsection{Construction of idempotent canonical relations}\label{CC}
We next recall how the condition $\mathcal{C}\circ\mathcal{C}=\mathcal{C}$ is managed.

Let $\mathcal{C}\subset T^*(X\times X)$ be a canonical relation, i.e. a conic Lagrangian Lagrangian sub-manifold of $T^*X\times T^*X^0$ ($^0$ means that we reverse the sign of the canonical symplectic form) such that $\mathcal{C}\circ\mathcal{C}=\mathcal{C}$. Then the positivity condition above implies that the projections $pr_1\mathcal{C}=\Sigma^+,pr_2\mathcal{C}=\Sigma^-$ are involutive complex sub-manifolds of $T*X$, with ${\rm Re\,}\Sigma^\pm=\Sigma$, $\Sigma^+\gg0$, $\Sigma^-\ll0$. $\gg0$ means that locally $\Sigma^+$ is defined by $n={\rm codim\,}\Sigma^+$ transversal equations $x_1=\dots x_n=0$ where the Poisson brackets $\{x_p,x_q\}$ vanish on $\Sigma^+$ and the matrix $(\frac1i\{x_p,\bar{x}_q\})$ is $\gg0$ on $\Sigma$.
The characteristic foliation of $\Sigma^\pm$ (tangent to the symplectic orthogonal $(T\Sigma^{\pm})^\perp$) is then tranversal to $\Sigma$, so it defines a projection $\Sigma^\pm\to\Sigma$, and we have $\mathcal{C}=\Sigma^+\times_\Sigma\Sigma^-$ (the linearized version of this is elementary).

\medskip There is a standard way of constructing such pairs of involutive manifolds $\Sigma^\pm$, cf.\cite{lB97.3}\,: 

\begin{Prop}\label{prC} Let $\delta$ be a smooth function on $T^*\Sigma$ vanishing of order $2$ on $\Sigma$, and such that ${\rm Re\,}\delta(\xi)\gtrsim {\rm dist\,}(\xi,\Sigma)^2$. Then there exists a unique pair $\Sigma^\pm$ such that $\delta$ vanishes on $\Sigma^+$ and $\Sigma^-$, $\Sigma\gg0$ and $\Sigma¬-\ll0$. $\Sigma^\pm$ is the outgoing, resp. in-going manifolds of the complex hamiltonian vector field $\frac1iH_a$ out of $\Sigma$. Any pair $\Sigma^\pm$ can be generated in that manner. The set of such functions $\delta$, so as those producing a given pair $\Sigma^\pm$, is contractible (convex). \end{Prop}

E.g. in the case of the Szeg\"o kernel on the sphere, $\Sigma^\pm$ is the complex characteristic manifold of $\bar{\partial}_b$, resp. $\partial_b$, and $\mathcal{C}$ is the complex flow of $\bar{\partial}_{b,z}\times\partial_{b,w}$ out of ${\rm diag\,}\Sigma$ in $T^*X\times T^*X$. The standard choice for $\delta$ is the symbol of $\square_b$.

\medskip \noindent Proof of Prop.\ref{prC}: we start with a smooth homogeneous function $\delta$ on $T^*\Sigma$ vanishing of order $2$ on $\Sigma$ and such that ${\rm Re\,}\delta(u)\ge c\  {\rm dist\,}(u,\Sigma)^2$. Let $H_\delta$ be its hamiltonian vector field. Then the the linearization along $\Sigma$ $\frac1iH_\delta $ is orthogonal to $\Sigma$ (for the symplectic structure), and its transversal part has no real eigenroot, so $\frac1i(H_\delta)$ has two well defined outgoing and in-going manifolds $\Sigma^+,\Sigma^-$ from $\Sigma$; these are involutive and with this choice of signs $\Sigma^+$ is $\gg0$ and $\Sigma^-$ is $\ll0$, as in the case of the Szeg\"o kernel of the sphere.
{\small(The corresponding result of linear algebra is the following: let $E$ be a real vector space with real symplectic form $\omega$, $q(x,y)$ a complex symmetric bilinear form such that $\delta(x)=\frac12{\rm Re\,}q(x,x)\ge c\,\|x\|^2$, $A$ the linear operator such that $q(x,y)=-\omega(Ax,y)$ (corresponding to the hamiltonian field of $\delta$): $A$ is antisymmetric w.r. to both $\omega$ and $q$; it has no real eigenroot, since if for some complex vector $z$ we have $Az=\lambda z$ with $\lambda$ real  ($\bar{A}\bar{z}=\lambda\bar{z}$), we have $q(z,\bar{z})=-\omega(Az,\bar{z})=-\lambda\,\omega(z,\bar{z})$,  $\bar{q}(\bar{z},z)=-\omega(\bar{A}\bar{z},z)=+\lambda\omega(z,\bar{z})$ hence ${\rm Re\,}q(z,\bar{z})=0$ and $z=0$. The spectral spaces $E_\lambda,E_{\lambda'}$ are orthogonal except for $\lambda'=-\lambda$ so the $E^\pm=\sum_{\pm{\rm Im}\lambda>0}E_\lambda$ are Lagrangian, in duality by $\omega$; $E^+$ is $\gg0$ an $E^-$ is $\ll0$ (${\rm Re\,}\frac1i\omega(z,\bar z)$ is negative on $E^+$ so for the set of linear forms vanishing on $E^+$ we get the other sign).
}

\medskip Conversely if the pair $\Sigma^\pm$ is given, there always exists a generating function $\delta$ as above; for instance locally $\Sigma^+$ can be defined by $n$ smooth equations $x_p=0$ with $\{x_p,x_q\}=0$, $\frac1i\{x_p,\bar{x}_q\}=\delta_{pq}=0$ ($\delta_{pq}$ the Kronecker symbol); $\Lambda^-$ is then defined smooth equations $\bar{y}_p=0$ with $\bar{y}_p=\bar{x}_p+\sum a_{pj}x_j+O(|x|^2)$ (mod. functions vanishing of second order on $\Sigma$). Since $\Sigma^-$ is involutive, the matrix $\alpha=(\alpha_{kj})$ is symmetric. The condition $\Sigma^-\ll0$ i.e. $\frac1i\{y_p,\bar{y}_q\}\gg0$ then means that the matrix $I-\alpha^*\alpha$ is $\gg0$, i.e. $\|a\|<1$. Then we set $\delta=\sum x_p\bar{y}_p$: its real part is $\ge(1-\|\alpha\|^2)|x|^2-O(|x|^3)$. We get $\delta$ globally by patching local results using a partition of unity (constructing $\delta$ real analytic would require a little more work but works just as well).

\medskip The function $\delta$ generating the pair $\Lambda^\pm$ is of course not unique, but clearly it can be chosen depending smoothly on a parameter if $\Lambda^\pm$ does so.

%%%%%%%%%%%  4.2  %%%%%%%%%%%%%
\subsection{Construction of the symbol (leading term) of a projector}\label{sS}

The set of Fourier integral operators associated with $\mathcal{C}$ is an algebra (without unit). The symbol of such operators lives on $\mathcal{C}$; choosing a frame, i.e. a basic elliptic symbol, identifies the set of symbols of these operators with the set of symbol functions on $\mathcal{C}$ (there is no canonical frame). We identify $\mathcal{C}$ with $\Sigma^+\!\times_\Sigma\Sigma^-$ and denote $(\xi',\xi")$ the variable, $\xi'\in\Sigma^+,\xi"\in\Sigma^-$ with common projection $\xi_0\in\Sigma$.
If $\sigma$ denotes the principal symbol in some frame, we have\,:
$$
\sigma(A\circ B)\;(\xi',\xi") = J(\xi',\xi")\;\sigma(A)(\xi',\xi_0)\;\sigma(B)(\xi_0,\xi")
$$
where $J$ is a fixed elliptic symbol. Since the composition law is associative we have 
$$
J(\xi',\xi_0)=J(\xi_0,\xi")=J(\xi_0,\xi_0)
$$

The change of frames $\tilde{\sigma}(A)=J(\xi,\xi")J(\xi_0,\xi_0)^{-2}\,\sigma(A)$ then gives 
\begin{equation}
\tilde{\sigma}(A\circ B)=\tilde{\sigma}(A)(\xi',\xi_0)\;\tilde{\sigma}(B)(\xi_0,\xi").
\end{equation}
If the frame is chosen in that manner ($J=1$), we have $a=\sigma(A)=\sigma(A\circ A)$ iff
\begin{equation}
a(\xi',\xi")=a(\xi',\xi_0)\,a(\xi_0,\xi") \qquad({\rm this\  requires}\  a(\xi_0,\xi_0)=1)
\end{equation}
Again there are many idempotent symbols, but it is obviously possible to keep track of this in a deformation.

%%%%%%%%%%%  4.3  %%%%%%%%%%%%%
\subsection{Construction of a projector mod. smoothing operators}\label{S}
\begin{Lemm} Let $S_0$ be a Fourier integral operator with canonical relation $\mathcal{C}$, whose symbol $a$ is a projector as above. Then the closed algebra generated by $S_0$ (closed for the filtration by degrees, mod. $C^\infty$) contains a unique projector $S$ such that $\rm{deg\,}S-S_0<0$.
\end{Lemm}
Proof\,: the solution is
\begin{equation} \label{proj}
S = S_0 + (2S_0-1) \sum_1^\infty \frac{(2k-1)!}{k!(k-1)!}\  R^k,\quad{\rm with }\ R=S_0-S_0^2
\end{equation}
(the degree of $R$ is  $<0$); this is equivalent to $(2S-1)=(2S_0-1)[(1-2S_0)^2]^{-\frac12}$.

\bigskip\noindent{\bf Deformation :}  let  $S_0,S_1$ be two Toeplitz projectors,  $\mathcal{C}_0,\mathcal{C}_1$ their canonical relations. It is clear from \ref{CC} that there exists a smooth family $\mathcal{C}_t$ of canonical relations as above linking them. It is clear from \ref{sS} that there exists a smooth family $s_t$ of symbol projectors linking the symbols of $S_0,S_1$ and a smooth family of Fourier integral operators $S'_t$ with symbol $s_t$ (not yet projectors). Formula \ref{proj} finally gives a smooth family of projectors linking $S_0$ and $S_1$ mod. $C^\infty$.

\bigskip\noindent {\bf Remark 1.} There in an analogue of Proposition \ref{prC} for Toeplitz projectors, that we do not use: the set of pseudo-differential operators such that $PS\sim0$, resp. $SQ\sim0$ is a positive left, resp. right, ideal $\mathcal{I}^\pm$, whose complex characteristic manifold is $\Sigma^\pm$, and $S$ mod. smoothing operators is uniquely determined by the fact that it is an elliptic Fourier integral projector, and $\mathcal{I}^+S\sim S\,\mathcal{I}^-\sim0$. There is no economic substitute for the function $\delta$ in that proposition;  $\square_b$ does play this role for the Szeg\"o projector ($\square_b S=S \square_b =0$) but in general the fact that an operator $P$ kills $S$, i.e. $P\in \mathcal{I}^-\mathcal{I}^+$, cannot be read on the principal symbol alone.

\bigskip\noindent {\bf Remark 2.} The argument above gives a projector mod. $C^\infty$; it does not require that $X$ be compact. Compacity is required to define the trace, and also to construct a true projector out of the approximate one, which is the useful thing for analysis. This construction is indicated in loc. cit. and we do not repeat it since it is not needed here. Modifying $S$ by a smooth projector of finite rank will of course not change its singularity, so there always remains an index ambiguity. 

If $X$ is of real dimension $3$ ($n=2$), and equipped with a pseudo-convex CR structure, the singularity of the Szeg\"o kernel is still well defined. The construction of \cite{BG81} or \cite{lB97.3} gives a Toeplitz projector with the right singularity, but this is not the Szeg\"o projector if $\bar\partial_b$ is not well behaved. The true Szeg\"o projector (i.e. the orthogonal projector on the space of holomorphic functions) is a Toeplitz projector essentially only if $X$ is holomorphically embeddable.

%%%%%%%%%%%  5  %%%%%%%%%%%%%
\section{Case of the 3-sphere}\label{sph}

In this section we examine the case where the base space is the oriented 3-sphere,  $X=S^3$\,:

\begin{Prop} 
For any contact form $\lambda$ on the 3-sphere, we have $L(\lambda)=0$
\end{Prop}
Proof: we identify the $3$-sphere with the unit quaternion sphere $SU_2$, set of all quaternions
$$
x=x_0+Ix_1+Jx_2+Kx_3 \quad\rm{with}\quad x\tilde{x}=x_0^2+x_1^2+x_2^2+x_3^2=1
$$
$\tilde{x}=x_0-Ix_1-Jx_2-Kx_3$ denotes the conjugate quaternion; to avoid confusions with complex conjugation, which is also needed. The complex coordinates are
$$
z_1=x_0+ix_1=r_1e^{i\theta_1},\quad z_2=x_2+ix_3=r_2e^{i\theta_2}
$$

The standard contact form is
\begin{equation}\label{standard}
\lambda_{st}=\langle dx,Ix\rangle = -i\, dz.\bar{z} =(x_0dx_1-x_1dx_0) + (x_2dx_3-x_3dx_2)= r_1^2d\theta_1+ r_2^2d\theta_2.
\end{equation}
We have $\lambda_{st}\,d\lambda_{st}=2\, dvol$, the canonical volume of the sphere; the conjugate form, giving  the opposite orientation, is 
$$
\tilde{\lambda}_{st}  = - r_1^2d\theta_1 +  r_2^2d\theta_2
$$

We make use of the following facts

\medskip \noindent{\bf 1.} Contact forms, just as three-manifolds, can be glued together: if $(X_1,\lambda_1),(X_2,\lambda_2)$ are two oriented contact manifolds, they can be glued together to give $(X_1\#X_2,\lambda_1\#\lambda_2)$ in the following manner: we choose base points $p_1,p_2$ and local coordinates so that a neighborhood of $p_1$, resp. $p_2$, is identified with the spherical cap $x_0<\frac12$ ($p_1$ the south pole $x=-1$), resp. $x_0>-\frac12$ ($p_2$ the north pole $x=1$), in the three sphere. We can deform $\lambda_1$, resp. $\lambda_2$, so that they coincide with the standard form on these coronas. Gluing is then obvious: $X_1\#X_2$ is obtained by gluing $X_1,X_2$ along the corona $-\frac12<x_0<\frac12$ using the identity map, and the contact forms patch together. Singularities of Szeg\"o or Toeplitz kernels can be glued in the same manner (remember that they live on the diagonal; they must be deformed so as to coincide with the standard Szeg\"o kernel on the two glued hemispheres). Since the logarithmic term of the standard Szeg\"o kernel vanishes, we get
\begin{equation} \label{add}
L(\lambda_1\#\lambda_2)=L(\lambda_1)+L(\lambda_2)
\end{equation}

\medskip \noindent{\bf 2.} Note that in real dimension $3$ a Toeplitz projector can always be chosen so as to correspond to a formally integrable boundary complex structure: the manifold $\Sigma^+$ is of codimension $1$ so the Frobenius integrability condition is empty, and its equation can always be chosen linear w.r. to $\xi$. So the Toeplitz projector is a "Szeg\"o projector  mod. $C^\infty$". This can be useful because the singularity of the (formal) Szeg\"o kernel is local and computable using Fourier integral operator calculus. Formally integrable tangent complex structures in dimension $3$ are usually not embeddable (in particular they never are when the corresponding contact form is over-twisted), so the Toeplitz projector analysis remains unavoidable. In higher dimension ($n\ge3$) tangent complex structures on a compact manifold are always embeddable (cf. \cite{lB75.2}), but contact structures do not necessarily contain any complex structure: what is gained on one side is lost on the other.

\medskip \noindent{\bf 3.} The cotangent bundle of $S^3$ is trivial. We identify it with $E\times S^3$, $V\simeq \mathbb{R}^3$ the space of pure imaginary quaternions, by $(x,u)\mapsto(x,ux)$, so that the standard contact form corresponds to the constant map $u(x)=I$. The map which to the homotopy class of a continuous function $a : S^3\to S^3=SU_2$ assigns the homotopy class of the vector field $u(x)=a(x)Ia(x)^{-1}x$, or of the map $x\mapsto u(x)=a(x)Ia(x)^{-1}\in S^2$, is the canonical isomorphism $\pi_3(S^3)=\mathbb{Z}\to \pi_3(S^2)$.

\medskip It follows from Eliashberg's analysis \cite{yE93,yE96,yE98}, that any homotopy class of such maps contains an over-twisted oriented contact form with the right orientation, and this is unique up to isomorphism (deformation, preserving the right orientation). The standard form corresponds to the constant map $a=1$ (the zero element of $\pi_3(S^2)$), so the trace integral defines a group homomorphism $L : \pi_3(S^2)=\mathbb{Z}\to\mathbb{R}$ ($L(\lambda)$ is real because the Toeplitz projector can be chosen self-adjoint).

\bigskip With the notations above, and with $(\frac{r_2}{r_1})^2=\tan\phi$ ($0\le\phi\le\frac\pi2$), we set 
\begin{equation}\label{ovtw}
\lambda_n=\cos(2n+1)\phi\, d\theta_1 + \sin(2n+1)\phi \,d\theta_2. 
\end{equation}
This is a contact form ($a(\phi)\,d\theta_1+b(\phi)\,d\theta_2$ is a contact form with the right orientation iff $ab'-ba'>0$, cf. \cite{dB82}); the coefficient of $d\theta_2$, resp. $d\theta_1$, must vanish for $\phi=0$, resp. $\phi=\frac\pi2$, so that $\lambda_n$ is the pull-back of a contact form on the 3-sphere. We have $\lambda_0=(\cos\phi+\sin\phi)\lambda_{st}$. 

\medskip As a non vanishing form, $\lambda_n$ is homotopic to $\cos\phi\,d\theta_1\pm\sin\phi\,d\theta_2$, with $\pm=\sin(2n+1)\frac\pi2=(-1)^n$: one can first deform linearly to $\lambda_n+\chi\ d\phi$ with $\chi$ a positive function vanishing near the ends, then still linearly, to $\chi d\phi+\cos\pi\,d\theta_1 \pm\sin\phi\, d\theta_2$, finally back to $\cos\phi\,d\theta_1\pm\sin\phi\,d\theta_2$ (in the deformation one must take care that the coefficient of $d\theta_1$, resp. $d\theta_2$, must vanish for $\phi=\frac\pi2$, resp. $0$; the sign at one end does not matter so long as the other coefficient is $\neq0$). Thus the homotopy class of $\lambda_n$ as a non-vanishing form is the trivial class for $n$ even, and $-1$, the class of $\tilde{\lambda}_{st}$ for $n$ odd.

\medskip Let us set $\psi=(2n+1)\phi$ ($0\le\psi\le(2n+1)\frac\pi2$), so $\lambda_n=\cos\psi d\theta_1+\sin\psi d\theta_2$. Because it can easily be deformed (cf. cor.\ref{soft}), we can construct our Toeplitz (or Szeg\"o) projector singularity so that it is the standard Szeg\"o projector for $0\le\psi<\frac\pi4$ or its pull-back by a suitable change of variables as below for $n\pi+\frac\pi4\le\psi\le(2n+1)\frac\pi2$ (for $n$ odd we get the complex conjugate complex structure corresponding to $-\lambda_{st}$), and so that is is invariant under the map $(\psi,\theta_1,\theta_2)\mapsto(\psi+\frac\pi2,-\theta_2,\theta_1)$ in between (as is $\lambda_n$). Then the logarithmic term vanishes near the ends and repeats itself periodically $2n$ times so we get\,:
\begin{equation}
L(\lambda_n)=nL(\lambda_1)
\end{equation}

Now it follows from the analysis of Eliashberg (cf. \eqref{add}) that we have $L(\lambda_n)=0$ if $n$ is even, so $L(\lambda_1)=0$, and finally $L(\lambda)=0$ for any $\lambda$ because the class of $\lambda_1$ is the negative generator of $\pi_3(S^2)$. 

\bigskip This ends the proof, and shows that the logarithmic trace vanishes for all contact forms on the three sphere.

%%%%%%%%%%%  bibli  %%%%%%%%%%%%%

\bigskip{\sc Universit\'e Pierre et Marie Curie, Analyse Alg\'ebrique, Institut de Math. de Jussieu,
Case 82 - 4, place Jussieu, 75252 Paris Cedex 05, France}

{\it E-pmail address:  boutet@math.Jussieu.fr}

\end{document}